\documentclass[12pt,a4paper]{article}

\usepackage[ruled,vlined]{algorithm2e}
\usepackage{algorithmic}

\usepackage{amssymb,amsmath}
\usepackage{latexsym}
\usepackage{color}
\usepackage[dvips]{graphicx}

\newtheorem{theorem}{Theorem}
\newtheorem{proposition}{Proposition}
\newtheorem{lemma}{Lemma}
\newtheorem{corollary}{Corollary}
\newtheorem{remark}{Remark}
\newtheorem{ex}{Example}
\newtheorem{definition}{Definition}
\newtheorem{conj}{Conjecture}

\newcommand{\re}{{\mathbb R}}

\newcommand{\n}{{\mathbb N}}

\newcommand{\cA}{{\cal{A}}}

\newcommand{\cI}{{\cal{I}}}

\newcommand{\cP}{{\cal{P}}}
\newcommand{\cM}{{\cal{M}}}

\newcommand{\be}{{\mathbf{e}}}
\newcommand{\bh}{{\mathbf{h}}}
\newcommand{\bep}{{\mathbf{p}}}

\providecommand{\com[2]}{\begin{tt}[#1: #2]\end{tt}}

\title{Stability of linear switching systems and Markov-Bernstein inequalities
for exponents\footnote{A preliminary version of this work has been presented at the IEEE CDC 2012 \cite{cdc-continuous}.}}

\author{Vladimir~Yu. Protasov and Rapha\"el~M. Jungers 
\thanks{The research of V.P. is supported by the
RFBR grants No
13-01-00642 and 14-01-00332, and by the grant of Dynasty foundation.  R. J. is an F.R.S.-FNRS Research Associate. His work is supported by the Communaut\'e fran\c caise de Belgique - Actions de Recherche
Concert\'ees, and by the Belgian Programme on Interuniversity Attraction Poles initiated
by the Belgian Federal Science Policy Office. }
\thanks{V. P. is with the Department of Mechanics and Mathematics, Moscow State University,
Vorobyovy Gory, 119992, Moscow, Russia. {E-mail: \tt\small
v-protassov@yandex.ru.}}%
\thanks{R. J. is with the ICTEAM Institute,
Universit\'e catholique de Louvain, 4 avenue Georges Lemaitre,
B-1348 Louvain-la-Neuve, Belgium. {E-mail: \tt\small
raphael.jungers@uclouvain.be.}}%
}
\date{}

\begin{document}

\maketitle

\begin{abstract}
We analyse the problem of stability of a continuous time linear switching system (LSS)
versus the stability of its Euler discretization. It is well-known that
the existence of a positive $\tau$ for which the corresponding discrete time system
with step size $\tau$ is stable implies the stability of LSS. Our main goal is to obtain a converse
statement, that is, to estimate the discretization step size $\tau > 0$ up to a given accuracy $\varepsilon > 0$. This leads to a method  of deciding the stability of continuous time LSS with a  guaranteed accuracy. As the first step, we solve this  problem for matrices with real spectrum and conjecture that our method stays valid for the general case.
Our approach is based on Markov-Bernstein type inequalities for systems of exponents. We obtain
universal estimates for sharp constants in those inequalities.

Our work provides the first estimate of the computational cost of the stability problem for continuous-time LSS (though restricted to the real-spectrum case).

\medskip

\noindent \textbf{Keywords:} {\em linear switching system, discretization step, stability, 
Lyapunov exponent, joint spectral radius, Markov-Bernstein inequality, Chebyshev system, Laguerre weight}
\smallskip

\begin{flushright}
\noindent  \textbf{AMS 2010} {\em subject
classification: }  37N35, 15A60, 34D08, 26D10, 41A50 
\end{flushright}

\end{abstract}

\begin{center}
\textbf{I. Introduction}
\end{center}
\smallskip

Linear switching systems (LSS) have been at the center of a great attention in the literature for the past years.  In continuous time, these are systems described by the following linear ODE:
\begin{equation}\label{main}
\left\{
\begin{array}{l}
\dot x (t) \ = \ A(t)\, x(t)\,  ;\\
x(0) \, = \, x_0\, ;\\
A(t) \in \cA \, ,\ t \ge 0,
\end{array}
\right.
\end{equation}
where $A(\cdot)$ is an arbitrary Lebesgue measurable function from $\re_+$ to
a given compact set of $d\times d$ matrices~$\cA$.
 On top of the theoretical challenges {they} offer, they appear in many applications such as viral disease treatments optimization~\cite{hmcb10}, multihop control networks \cite{PappasTAC2011}, multi-agent and consensus systems \cite{AC,BHT}, etc. 
They also provide a theoretical and computational framework for analyzing more complex systems, like general hybrid systems, systems with quantized signals, or event-triggered control schemes, etc. (see~\cite{Lib} for a general survey).

System  (\ref{main}) is {\em stable} if
 $x(t) \to 0$ as $t \to +\infty$, for every initial condition ${x_0 \in \re^d}$
  and any measurable function $A(\cdot): \re_+ \to \cA$.
   The problem of deciding whether a switching system is stable has been studied in many papers
    (see bibliography in~\cite{Lib}).
  The stability analysis amounts to compute the so-called {\em Lyapunov exponent} (also called {\em worst-case Lyapunov exponent}) of the corresponding set of matrices:
	\begin{definition} The \emph{worst-case Lyapunov exponent} of System \ref{main} is defined as the infimum of the set of values $\alpha$ such that
	for all trajectories of~(\ref{main}), $$\exists C:\, \|x(t)\| \, \le \, C e^{\, \alpha t}.$$

\end{definition} Obviously, if $\sigma < 0$, then the system is stable. The converse is also true~\cite{O,MP}, and so, the stability of  LSS
defined by the set of matrices $\cA$ is equivalent to the
  condition  $\sigma(\cA) < 0$. The stability analysis of continuous time LSS is an NP-hard problem \cite{gurvits-olchevski}, but in fact, no method is known that approximates the Lyapunov exponent, even in exponential time, up to a guaranteed accuracy. The methods proposed in the literature are, to the best of our knowledge, only sufficient conditions for stability, and this is even more surprising in view of the fact that the equivalent question for discrete-time switching systems has found a positive answer since quite a long time (see \cite[Section 2.3]{jungers_lncis} for a survey on that question).
  It may happen, however, that neither of those conditions are satisfied, although the system is very stable
  ($\sigma(\cA) << 0$).  The problem of approximating $\sigma (\cA)$ with a prescribed accuracy seems to be very hard. In this paper, we make the first step
towards its solution and tackle the special case where all the matrices from~$\cA$ have real spectra.  In this case, we are able to compute the step size $\tau$ such that the stability of the discretized switching system (i.e., defined by matrices $I+\tau \cA=\{I+\tau A:A\in \cA\}$) is equivalent, up to a given precision $\varepsilon > 0$, to the stability of
the original system~(\ref{main}). Thus, the algorithm reduces the stability issue
of an LSS to a discrete time LSS, for which efficient methods for deciding stability are known.
Surprisingly, the guaranteed step size $\tau$ is not too small:
in decays linearly with $\varepsilon$ and quadratically with the dimension~$d$. This makes our algorithm applicable in practice as demonstrated in a numerical example in Section~6.  We believe that the same result holds for general systems as well, although our method of proof only works for matrices with real spectra.

Our technique to analyse the discretization of LSS relies on Markov-Bernstein inequalities
for exponential polynomials. These are inequalities of the type $\|p^{(k)}\|\le C \|p\|$, where
$p = \sum_{j=1}^d p_je^{-h_jt}$ is a polynomial of exponents with given positive parameters $\{h_i\}_{i=1}^{k}$.
For every $d$ and $k$, we derive uniform upper bounds for the constant $C$, over all polynomials with parameters
$h_i \in (0, 1]$ (Theorem~\ref{th3}). Combining this with recent results of Sklyarov~\cite{S}
on Markov-Bernstein inequalities with the Laguerre weight, we estimate that constant in terms of the dimension~$d$
and the order of the derivative~$k$  (Theorem~\ref{th4}). In Sections~4 and~5 we apply these results to
obtain uniform lower bounds for the step size~$\tau$ in case of matrices with real eigenvalues.
This is formulated in Theorems~\ref{th5} and~\ref{th10}. We conjecture that the same holds for general matrices.

\smallskip

\begin{center}
\textbf{2. Statement of the problem}
\end{center}
\smallskip

We consider a linear switching system (LSS) of the form (\ref{main}).
If $A(t)$ is a measurable function taking values in a given compact set of $d\times d$ matrices~$\cA$, then
the solution of~(\ref{main}), a univariate absolute continuous vector-function~$x(t)$  taking values in~$\re^d$, is called a {\em trajectory} of the system. So, the system is stable if all its trajectories converge to zero as $t \to \infty$.
A well-known  necessary condition for stability is that each matrix $A \in \cA$ is Hurwitz, i.e., all its
eigenvalues have strictly negative real parts. Indeed, if it is not the  case, then the trajectory $x(t) = e^{t A}x_0$
corresponding to the constant function $A(t) \equiv A$ does not tend to zero as $t \to \infty$. In the sequel we
assume this condition is satisfied.
\smallskip

We use the following notation: $\re_+^d = \{x \in \re^d\, , \, x \ge 0\}, \
\re_{++}^d = \{x \in \re^d \, , \,  x > 0\}\, , \, I$ is the identity matrix, ${\rm sp}(A)$ is the spectrum
of the matrix~$A$ (the set of eigenvalues counted with multiplicities),
$\rho(A)$ is the spectral radius (the largest modulus of eigenvalues), $\be = (1, \ldots , 1) \in \re^d$ is the vector of ones.
\smallskip

One method for deciding the stability of LSS elaborated in the literature is to pass
to the corresponding {\em discrete time system} with a step size~$\tau > 0$:
\begin{equation}\label{main-d}
X(k+1)\ = \  X(k) \, + \, \tau \, A(k)\, X(k)\, , \quad A(k) \in \cA\, , \ \forall k\in \n\, , \  X(0) = X_0\in \re^d\, .
\end{equation}
This system is obtained from~(\ref{main}) by replacing~$x(k\tau)$ with $X(k)$ and  the derivative~$\dot x (k \tau)$
with the divided difference~$\frac{X(k+1) \, - \, X(k)}{\tau}$. 
The stability of the discrete system, which
means $X(k) \to 0$ for every $X_0$ and for every sequence of matrices $A(k) \in \cA$,
is equivalent to the condition~$\rho(I +\tau \cA)< 1$, where $\rho(I +\tau \cA)$ is the {\em joint spectral radius}
of the family~$I + \tau \cA \, = \, \{I + \tau A \ | \  A \in \cA\}$.
\begin{definition}\label{d20}
The  joint spectral radius (JSR) of a set of matrices $\cM$ is
\begin{equation}\label{eq:JSR}
\rho(\cM)\quad =\quad \lim_{k \to \infty} \
\max_{A_1,\dots,A_k\in \cM}\ \|A_1\cdots A_k\|^{1/k},
\end{equation}
where $\|\cdot\|$ is an arbitrary matrix norm.
\end{definition}
Thus, the stability problem for a discrete time LSS
is reduced to the problem of computing the corresponding joint spectral radius $\rho(I + \tau A)$. Although this
problem is known to be hard~\cite{BT}, some efficient algorithms have been constructed in the last years
for approximate computation of JSR.  See \cite{jungers_lncis} for a general survey on the topic. In particular, approximation algorithms with guaranteed accuracy are available. Moreover,
a recent line of work~\cite{GP,JCG} even presents algorithms
 that find the exact value of JSR (in the form of a root of some polynomial) for a vast majority of finite matrix families of dimensions up to $d = 20$ and higher.

Our goal is to estimate the parameter~$\tau$ to be able to infer the stability
of a continuous time system from its discretization.
The following crucial fact was proved in the eighties:
\smallskip

\noindent \textbf{Theorem A.}~\cite{O, MP}. {\em If there exists~$\tau > 0$ for which the discrete time system~(\ref{main-d}) is stable, then it is stable for all smaller positive~$\tau$, and the continuous system~(\ref{main})
is stable.}
\smallskip

It is natural that the stability of the continuous time LSS
follows from the stability of its discrete approximation, provided $\tau > 0$ is small enough.
What is more surprising is that by Theorem~A, this is true for a particular
step size $\tau$, not necessarily a small one. If there exists $\tau > 0$ such that~$\rho(I + \tau \cA) < 1$, then the continuous time LSS is stable. The converse is also true:
\smallskip

\noindent \textbf{Theorem B.}~\cite{O, MP}. {\em  If the continuous time LSS is stable,
then there is~$\tau > 0$ such that the corresponding continuous system~(\ref{main}) is stable, i.e., ~$\rho(I + \tau \cA) < 1$.}
\smallskip

The practical implementation of Theorems A and B
may be hard, because they do not specify the step size~$\tau$.  If we take, say, $\tau = 10^{-3}$
and get~$\rho(I + \tau \cA) > 1$, then no conclusion can be drawn on the stability of the continuous~LSS.
The inequality~${\rho(I + \tau \cA) < 1}$ can still be true for some smaller~$\tau$. The problem can now be formulated as follows:
\smallskip

{\em For which $\tau > 0$ does the inequality~$\rho(I + \tau \cA) \ge 1$ imply} that $\sigma \ge 0$, i.e., does imply instability~?

\smallskip

This problem, however, is not well-posed, because such a universal step size~$\tau$, for all compact sets of matrices~$\cA$,
does not exist. Indeed, if one multiplies all matrices from~$\cA$ by a large number $N$, then the
step size is replaced by $\tau/N$ which tends to zero as $N \to \infty$. Hence, the family $\cA$ has to be normalized.
Taking as a normalization factor the largest spectral radius of matrices from $\cA$ we obtain the second version of the
problem:
\smallskip

{\em For which $\tau (r) > 0$ does the following hold}:
if $\max_{A \in \cA} \rho(A) \le r$, then the inequality~${\rho(I + \tau \cA) \ge 1}$ implies that $\sigma \ge 0$~?
\smallskip

 For this problem, there is no positive answer either. It is easy to construct families of matrices with $r=1$ and
 with arbitrary small step size $\tau$.  Moreover, no algorithm is known, to the best of our knowledge, which allows to answer the question \lq $\rho<1?$ \rq in finite time \cite[Section 4.1]{jungers_lncis}.  A well-posed formulation should require that the stability of the LSS is only determined up to a prescribed accuracy~$\varepsilon > 0$.
\smallskip

\noindent \textbf{Problem 1.} {\em Given $r > 0$ and $\varepsilon >0,$
Find the maximal step size $\tau = \tau(r, \varepsilon)$ for which the following holds:
if $\cA$ is an arbitrary compact family of matrices such that $\max_{A \in \cA} \rho(A) \le r,$ then the inequality~$\rho(I + \tau \cA) \ge 1$ implies that $\sigma(\cA) \ge -\varepsilon .$}

\smallskip

\begin{definition}\label{d10}
For  given parameters $r,\varepsilon > 0$, the value
$K(r, \varepsilon)$ is the largest positive number such that, for every~$\tau < K(r, \varepsilon)$ and for any set of matrices $\cA$ such that $ \max\, \{\rho(A) \ | \ A\in\cA\} \, \le \, r$,
the inequality~$\, \rho(I + \tau \cA)\ge 1\, $ implies $\, \sigma (\cA) \, \ge \, - \varepsilon.$
\end{definition}

  Thus, for every~$\tau < K(r, \varepsilon)$ we have: if~$\, \rho(I + \tau \cA)< 1\, $, then $\sigma (\cA) < 0$
(the LSS $\dot x = Ax\, , \ A(t) \in \cA$ is stable); otherwise, if~$\, \rho(I + \tau \cA)\ge 1\, $, then $\sigma(\cA) \ge -\varepsilon$ (the LSS $\dot x = (A+\varepsilon I)x\, , \ A(t) \in \cA$ is unstable)\footnote{Note that, as said above, the very question $\rho(I + \tau \cA)< 1$ cannot be solved in finite time.  However, we will show in Section 6 that an answer to a relaxed version of this question is sufficient for our purposes.}. To solve Problem~1 one needs to find a computable lower bound for $K(r, \varepsilon)$. This would allow us not only to determine the stability/instability of an LSS but
to also evaluate its Lyapunov exponent $\sigma (\cA)$ with a given precision.
Indeed, since for an arbitrary number $a \in \re$, we have $\sigma (\cA - a I) = \sigma (\cA) - a$ (see, for instance,~\cite{MP}), it follows that we can decide between two cases: $\sigma(\cA) < a$ and $\sigma (\cA) \ge a - \varepsilon$, just by computing the joint spectral radius $\rho\bigl(I + \tau (A - a I)\bigr)$ for some $\tau < K(r, \varepsilon)$.
Choosing suitable sequence of numbers $a$, one can compute $\sigma (\cA)$ by double division.

Thus, do decide the stability and to compute the Lyapunov exponent $\sigma (\cA)$ by means of discretization
one has to know a lower bound for $K(r, \varepsilon)$. Moreover, this quantity should not be too small,
otherwise computing $\rho(I + \tau A)$ becomes hard, because all matrices are close to the identity matrix.
An obvious way to estimate $K(r, \varepsilon)$ is by using  the local Lipschitz continuity of the joint spectral
radius. However, all bounds for the local Lipschitz constant of JSR available in the literature~(see~\cite{P, Koz})
grow exponentially in $d$ and may be very large if some matrices from $\cA$ have two close eigenvectors.
That is why this technique leads to estimates for $K(r, \varepsilon)$ that are too small and not practical.

In this paper, we suggest a different approach based on a geometrical analysis of the Lyapunov norm of the family~$\cA$.
This leads to lower bounds for $K(r, \varepsilon)$ in terms of the sharp constants in the Markov-Bernstein
inequalities for exponents. In the next section we give a short overview on Markov-Bernstein inequalities
and prove universal upper bounds for the constant in these inequalities
(Theorems~\ref{th3} and \ref{th4}). Next, in Section~4, we formulate the main results of the paper --
Theorems~\ref{th5} and \ref{th10} that estimate the discretization step size $\tau$ of LSS in terms of the bounds from Section~3. The proofs and further estimates  are given in Section~5, and Section~6 presents in details our algorithm and a numerical example.

\smallskip

\begin{center}
\textbf{3. Markov-Bernstein inequalities for exponents}
\end{center}
\smallskip

Let $\Phi = \{\varphi_i (\cdot)\}_{i=0}^n$ be a collection of continuous functions on an interval
$\cI \subset \re$ which is either a segment, or half-line, or the whole real line.
We denote by $\cP = \cP(\Phi)$ the linear subspace of $C(\cI)$
spanned by elements of~$\Phi,$ which we refer to as the space of \lq polynomials\rq{}  of this system.
For a given $k \in \n$ we consider the value
$$
C_k \ = \ C_k(\Phi) \ = \ \max\, \bigl\{ \|p^{(k)}\|_{\cI} \quad \bigl| \quad  p \in \cP, \, \|p\|_{\cI} \le 1\bigr\},
 $$
where $||f||_{\cI}\triangleq \max_{x \in {\cI}}{|f(x)|}.$
 Thus, $C_k$ is the maximal absolute value of the $k$th derivative
of functions from the unit ball of~$\cP$. From the compactness argument it follows that
the maximum is always attained. The corresponding inequalities are called {\em Markov-Bernstein inequalities} for the system $\Phi$:
\begin{equation}\label{markov1}
\forall p \in \cP, \quad \bigl\|p^{(k)}\bigr\|_{\cI} \ \le  \  C_k \, \bigl\|  p \bigr\|_{\cI}.
\end{equation}
This type of inequalities and its generalizations (weighted inequalities, Kolmogorov's type inequalities, etc.)
between functions and their derivatives, are popular topic in approximation theory, real analysis,
and optimization. It has been studied in the literature in great detail (see~\cite{BE0, BE3, Nat, MMR, Tikh} and the references therein).
The best known of them are the classical Bernstein and Markov inequalities. The Bernstein inequality
states that for trigonometric polynomials of degree at most $n$
(i.e., in the case $\varphi_m(t) = e^{i m t}, m = 0, \ldots , n$) we have
$\|p'\|_{[0, 2\pi]} \, \le \, n \, \|p\|_{[0, 2\pi]} $. The constant $n$ is attained
for $p(t) \, = \, c\,  e^{i n t}$ only. The Markov inequality is for algebraic polynomials
of degree at most $n$   (the case $\varphi_m(t) = x^m, m = 0, \ldots , n$). It states that
$\|p'\|_{[-1, 1]} \, \le \, n^2 \, \|p\|_{[-1,1]}$, and the constant is attained
for the corresponding Chebyshev polynomial of degree~$n$. The sharp constants and extremal polynomials in the Markov and Bernstein inequalities are known for derivatives of all orders $k$ and their properties are well studied.
We will need such inequalities for the system of
real exponents $\varphi_i(t) = e^{-h_i t}, \, h_m > 0$, on the half-line~$\re_+$, for which less in known.
We start with introducing some notation.

In the sequel, it will be more convenient to enumerate functions from $1$ to $d$.
We consider a vector $\bh = (h_1, \ldots , h_d) \in \re^d_+$ such that
$0< h_1 \le \cdots \le h_d$ and the corresponding system of
exponents $\Phi_{\bh} = \{e^{-h_i t}\}_{i=1}^d$. The space of polynomials of this
system on the half-line $\re_+$ will be denoted by $\cP_{\bh}$. This is a $d$-dimensional subspace of the space
$C_0(\re_+)$ of functions continuous on $\re_+$ and converging to zero as $t \to +\infty$.
We allow some of the numbers~$h_i$
to coincide, in which case the corresponding exponents
 are multiplied by powers of $t$:  if $h_i < h_{i+1} = \cdots = h_{i+m}  < h_{m+1}$, then
the exponent $h_{i+1}$ has multiplicity $m$ and the functions $e^{-h_{i+1}t}, e^{-h_{i+2}t} \ldots , e^{-h_{i+m} t}$
are replaced by $e^{-h_{i+1}t}, te^{-h_{i+1} t}, \ldots , t^{m-1}e^{-h_{i+1} t}$ respectively.
The map $\bh \mapsto \Phi_{\bh}$ is thus well-defined and continuous~\cite{K, KN}.

It is well known that for every $\bh$ the system of exponents $\Phi_{\bh}$
 is a Chebyshev system on~$\re_+$, i.e., every nontrivial polynomial from $\cP_{\bh}$
has at most $d-1$ zeros~(see, for instance,~\cite{KN,KS, Dz}).
As a Chebyshev system, it has the following properties. For any set of $d$ distinct points~$t_i \in \re_+$ and for any set of numbers
$c_1, \ldots , c_d$,  there is a unique polynomial~$p \in \cP_{\bh}$ such that $p(t_i) = c_i\, , \, i = 1, \ldots , d$.
 By Haar's theorem~\cite{KN},
for every continuous function  $f\in C(\re_+)$,  there is  a unique element $p \in \cP_{\bh}$  of best approximation,
for which the value $\|f - p\|_{\re_+}$ attains its minimum on the set $\cP_{\bh}$. By Karlin's theorem (the ``snake theorem'', \cite{K}),
the difference $f - p$ is either  identically zero, or it possesses $d$ points of Chebyshev alternance, where $f-p$ takes values equal by module with alternating signs. There exists a unique polynomial $T = T_{\bh} \in \cP_{\bh}$ and a unique system of points
$0= \nu_1 < \nu_2 < \cdots < \nu_d < \infty $ such that $\|T\|_{\re_+} = 1$ and $T(\nu_k) = (-1)^k, \, k = 1, \ldots , d$. We call $T_{\bh}$ the {\em $\bh$-Chebyshev polynomial}. This is a polynomial from $\cP_{\bh}$ with the smallest deviation from zero among all polynomials with a given leading coefficient
(i.e., coefficient for $t^{m-1}e^{-h_dt}$, where $m$ is the multiplicity of~$h_d$). This polynomial enjoys many extremal properties on the unit ball $\{p \in \cP_{\bh} \, , \  \|p\|_{\re_+} \le 1\}$.
In particular, for any $k \in \n$, the Chebyshev polynomial $T_{\bh}$ is a unique (up to the sign) solution of the problem
\begin{equation}\label{markov2}
\left\{
\begin{array}{l}
\max \ \bigl\|p^{\, (k)}\bigr\|_{\re_+} \\
\mbox{s.t. }\ p \in \cP_{\bh}\, , \ \|  p \|_{\re_+} \le 1\, .
\end{array}
\right.
\end{equation}
  The optimal value of this optimization problem will be denoted by $M_k(\bh)$.
  It is attained for the Chebyshev polynomial $T_{\bh}$ at $t=0$
  (see~\cite{KS} for the proofs).
  Thus, $M_k(\bh) = (-1)^{k+1}T_{\bh}^{(k)}(0)$.
  This value is the best possible constant in the {\em Markov-Bernstein inequality for exponents}:
\begin{equation}\label{markov3}
\bigl\|p^{\, (k)}\bigr\|_{\re_+} \ \le  \  M_k(\bh) \, \bigl\|  p \bigr\|_{\re_+} \, , \quad p \in \cP_{\bh}\, .
\end{equation}
Apart from a few special cases (for instance, when all $h_i$ are equal, or when they constitute an
arithmetic progression), the values $M_k(\bh)$ are not known. They, however, can be evaluated
numerically for each~$\bh$, by an approximate computation of the corresponding Chebyshev polynomial
 $T_{\bh}$ using the Remez algorithm~\cite{R1, Dz}. Lower and upper bounds for $M_1(\bh)$ in terms of the sum~$\sum_{i=1}^d h_i$ were obtained by Newman~\cite{New}, see also~\cite{BE1, BE2, BE3} for further generalizations.

 For our applications to linear switching systems,
 we  need a uniform estimate for the constants~$M_k(\bh)$ over the polytope (actually, simplex)
  ${\Delta_d = \bigl\{\bh \ge 0 \ \bigl| \  0< h_1 \le \ldots \le h_d \le 1 \bigr\}}$.
 We omit the index~$d$ if the dimension is specified.
 Thus, $\Delta$ consists of ordered positive vectors for which $h_d\le 1$. Since
 for any $\lambda > 0$, we have $M_k(\lambda \, \bh) = \lambda^k M_k(\bh)$, it suffices
 to estimate the constants $M_k(\bh)$ for $\bh \in \Delta$. In what follows, we assume
 $\bh \in \Delta$. Recall that  $\be\in \re^d$ is the vector of ones. The following theorem  establishes a sharp upper bound
 for this constant over all $\bh \in \Delta$.
 \begin{theorem}\label{th3}
The value $M_k(\bh)$ attains its maximum on the set $\Delta$ at a unique point $\bh = \be$.
\end{theorem}
This theorem is analogous to comparison theorems for systems of hyperbolic sines for~$k=1$~\cite{BE2}, but the method of proof is different.
We give the proof in Appendix, along with another comparison type result, Theorem~\ref{th100}, which is crucial in Section~5.

Let us denote $M_k(\be) = M_{k, d}$. Thus, the biggest possible constant $M_k(\bh)$
corresponds to the case when all~$h_i$ are maximal. Since $h_1 = 1$ has multiplicity $d$, the space $\cP_{\be}$
consists of exponential polynomials of the form~$p(t) = e^{-t}q(t)$, where
$q$ is an algebraic polynomial of degree at most $d-1$. The corresponding Chebyshev polynomial
$T_{\be}$ will be denoted by $S_d(t) = e^{-t}s_d(t)$, where $s_d$ is an algebraic polynomial.
Thus, $M_{k, d} = (-1)^{k+1} S^{(k)}_d(0)$.
\begin{corollary}\label{c5}
For every $k \in \n, \, \bh \in \Delta_d$ and an exponential polynomial $p \in \cP_{\bh}$, we have
$$
\bigl\| p^{(k)}\bigr\|_{\re_+} \ \le \ M_{\, k, d}\, \bigl\| p\bigr\|_{\re_+}\, .
$$
The equality is attained at a unique (up to normalization) polynomial $p = S_d$.
\end{corollary}
Asymptotically sharp upper bounds for $M_{k, d}$ have been derived in the literature.
It turns out that the algebraic polynomial  $s_d(t)$ is the solution
of the following extremal problem: among all algebraic polynomials $q(t)$ of degree at most~$d-1$ such that
$\|e^{-t} q(t)\|_{\re_+}\le 1,$ find the maximal value of~$\|q^{(k)}\|$.
It was show in~\cite{CLM, MN} that for all~$k \in \n$, a unique solution is given by the polynomial
$s_d$. Writing $L_{k, d}$ for the value of this problem, we
obtain the {\em Markov inequality for the Laguerre weight}:
\begin{equation}\label{markov4}
\bigl\|\, q^{(k)}\, \bigr\|_{\re_+} \quad \le  \quad  L_{k, d} \ \bigl\| \, e^{-t} q(t) \, \bigr\|_{\re_+} \, ,
\end{equation}
for every algebraic polynomial~$q$ of degree $\, \le \, d-1$.
In contrast to the classical Markov inequality, the norm of the polynomial is measured with the
{\em Laguerre weight} $e^{-t}$.  The extremal polynomial $s = s_{d-1}$ is called
the {\em Chebyshev polynomial with the Laguerre weight}.
It is characterized by existence of $d$ points of alternance
$0 = \nu_1 < \ldots < \nu_d$ such that $s(\nu_k) = s\, '(\nu_k) = (-1)^ke^{\, \nu_k}, \, k = 2, \ldots , d$.
We have $L_{k, d} = |s^{(k)}_d(0)|$. The first upper bound for this quantity was obtained in 1964 by
Szeg\"o~\cite{Sz}, who proved that $L_{1, d} \le Cd$, then this result was sharpened
in~\cite{Fre}, see also~\cite{CLM, MN}. Sklyarov in 2010 obtained a comprehensive solution to this problem:
\smallskip

\noindent \textbf{Theorem C}~\cite{S}. {\em
For every $k \le d-1$ we have
\begin{equation}\label{skl}
\frac{8^k(d-1)\, !\, k\, !}{(d-1-k)\, !  \, (2k)\, !} \Bigl(1 - \frac{k}{2(d-1)} \Bigr)\quad \le \quad L_{k, d}\quad \le \quad \frac{8^k(d-1)\, !\, k\, !}{(d-1-k)\, !  \, (2k)\, !}.
\end{equation}
}
\smallskip

This gives an upper bound for $L_{k, d}$ which is, moreover, asymptotically tight as $d \to \infty$ and $k$ is fixed.
On the other hand, as it was noted in~\cite{MN}, the value $M_{k, d}$ can be expressed with
the constants $L_{j, d} , \ j = 1, \ldots ,  k$,  as follows:
\begin{equation}\label{LM}
M_{k, d} \quad = \quad 1 \ + \ \sum_{j=1}^k \, \, {k \choose j}\ L_{j, d} \, .
\end{equation}
Combining this with~(\ref{skl})
we obtain after elementary simplifications:
\begin{equation}\label{LM1}
M_{k, d} \quad \le \quad  \sum_{j=0}^k \, \, \frac{8^j \ {d-1 \choose j} \ {k \choose j}}{{2j \choose j}} \,
\end{equation}
(all terms with $j \ge d$ are zeros;  ${n \choose 0} = 1$ for any $n \ge 0$). We did not succeed in any further
simplification of this expression. Combining with Theorem~\ref{th3}, we obtain:
\begin{theorem}\label{th4}
For every $k, d\in \n$ we have
$$
\max_{\|p\| \le 1, \ p  \in \cP_{\bh}, \ \bh \in \Delta_d}
\ \bigl\| \, p^{\, (k)}\, \bigr\| \quad \le \quad  \sum_{j=0}^k \, \, \frac{8^j {d-1 \choose j} {k \choose j}}{{2j \choose j}}
$$
(all the terms with $j\ge d$ are zeros).
This inequality is asymptotically tight as $d \to \infty$, the extremal polynomials are $p = S_d$.
\end{theorem}
In the next sections we need this result only for $k=2$, when~(\ref{LM1}) reads
\begin{equation}\label{k2}
M_{2, d} \ \le  \ \frac{16 \, d^{\, 2} \, - \, 24 \, d \, + \, 11}{3}\, .
\end{equation}
This upper bound is sharp only asymptotically, and actual
values of $M_{2, d}$ are smaller. They were listed in~\cite{S} for all $d = 2, \ldots , 20$.
In table~1 we write~$M_{\, 2, d}$ for~$d \le 10$.

\begin{table}[thb]
\begin{center}
\begin{tabular}{|c|c|}\hline
 $d$ & $M_{2, d}$ \\
 \hline
\rule{0pt}{9pt}\noindent
 $2$ &  8.182 \\
  \hline
 $3$ &  25.157 \\
  \hline
 $4$ &  52.587 \\
  \hline
 $5$ &  90.585\\
  \hline
 $6$ &  139.191 \\
  \hline
 $7$ & 198.420\\
  \hline
 $8$ &  268.283 \\
  \hline
 $9$ &  348.788\\
  \hline
 $10$ & 439.938 \\
 \hline
\end{tabular} \\[2mm]
\caption{\footnotesize{The values of~$M_{2, d}$ for $d = 2, \ldots , 10.$}}\label{tab1}
\end{center}
\end{table}
\smallskip
\newpage 

\begin{center}
\textbf{4. The main results}
\end{center}
\smallskip

In the special case when all matrices from $\cA$ have real eigenvalues, the results of previous section
generate lower bounds for the discretization parameter~$\tau$ of LSS. We formulate here the fundamental theorems; their proofs
along with other results, are given in the next section.

The following theorem allows to pick the discretization step size $\tau$ providing an accuracy~$\varepsilon,$ given the
dimension~$d$ and the maximal spectral radius $r$ of the matrices.

\begin{theorem}\label{th5}
If all matrices of~$\cA$ have real spectra, then
$$
 K(r, \varepsilon) \quad \ge \quad \frac{2\, \varepsilon}{M_{\, 2, d}\ r^2}\, .
$$
\end{theorem}
Invoking the upper bound~(\ref{k2}) for $M_{\, 2, d}$, we obtain
\begin{theorem}\label{th10}
If all matrices of~$\cA$ have real spectra, then
$$
\, K(r, \varepsilon) \ \ge \ \frac{6\, \varepsilon}{\bigl( 16\, d^2 \, - \, 24 d\, +\, 11\bigr)\, r^2}\, .
$$
\end{theorem}
\begin{corollary}\label{c10}
Suppose all matrices of~$\cA$ have real spectra; then if
the discrete time system with the step length
\begin{equation}\label{recipe}
\tau \, = \, \frac{6\, \varepsilon}{\bigl( 16\, d^2 \, - \, 24 d\, +\, 11\bigr)\, r^2}
\end{equation}
 is not stable, then~${\sigma (\cA) \ge - \varepsilon}$, i.e., the continuous time system
with the set of matrices $\cA + \varepsilon I$ is not stable.
\end{corollary}
\smallskip

We see that the lower bound~$K(r, \varepsilon)$ for the critical value of~$\tau$ is linear in~$\varepsilon$, which is natural, and decays with the dimension as~$d^{\, -2}$, which is much better than one could expect.

Certainly, the real spectrum assumption in our results is very restrictive. That is why we consider them as the first step towards the solution
of the problem of estimating the discretization parameter. Actually, we believe that Theorem~\ref{th10} holds for general matrices.
\begin{conj}\label{conj10}
Theorem~\ref{th10} and Corollary~\ref{c10} are true for the general case, without the real spectra assumption.
\end{conj}
In Section~5 (Remark~\ref{r30}) we discuss this conjecture.
\smallskip

Substituting in Theorem~\ref{th5} we get the lower bounds shown in table~\ref{tab2}.

\begin{table}[thb]
\begin{center}
\begin{tabular}{|c|l|}\hline
 $d$ & $\mbox{Lower bound for}~K(r, \varepsilon)$ \\
 \hline
\rule{0pt}{9pt}\noindent
 $2$ &  $\frac{1}{4.1 \, r^2}\ \varepsilon$ \\
 \hline
 $3$ &  $\frac{1}{12.6 \, r^2}\ \varepsilon$ \\
  \hline
 $4$ &  $\frac{1}{26.3 \, r^2}\ \varepsilon$ \\
  \hline
 $5$ &  $\frac{1}{45.3 \, r^2}\ \varepsilon$\\
  \hline
 $6$ &  $\frac{1}{69.6 \, r^2}\ \varepsilon$ \\
  \hline
 $7$ & $\frac{1}{99.3 \, r^2}\ \varepsilon$\\
  \hline
 $8$ &  $\frac{1}{134.2 \, r^2}\ \varepsilon$ \\
  \hline
 $9$ &  $\frac{1}{174.4 \, r^2}\ \varepsilon$\\
  \hline
 $10$ & $\frac{1}{220 \, r^2}\ \varepsilon$ \\
 \hline
\end{tabular} \\[2mm]
\caption{\footnotesize{Lower bounds for $\, K (r, \varepsilon)\, $ in Theorem \ref{th5}}}\label{tab2}
\end{center}
\end{table}
These bounds are better than those from Theorem~\ref{th10} because  they are based on the actual values of
 $M_{2, d}$ from table~\ref{tab1}, while Theorem~\ref{th10} uses the general upper bound~(\ref{k2}).

%
%
%

\begin{center}
\textbf{5. Individual estimates and proofs}
\end{center}
\smallskip

In this section we prove Theorem~\ref{th5} by first introducing the individual
maximal step size $\tau = K(\cA, \varepsilon)$ for a given family $\cA$ and then
obtaining the universal bound $K(r, \varepsilon)$ merely by taking infimum over all
matrix families $\cA$ with the largest spectral radius~$r$.
\begin{definition}\label{d30}
For  a given compact family of matrices $\cA$ and for $\varepsilon > 0$ let
\begin{equation}\label{Kind}
K(\cA, \varepsilon) \ = \
\left\{
\begin{array}{lcl}
+\infty & \ , \ & \sigma(\cA) \ge - \varepsilon\, , \\
\sup\, \bigl\{\tau > 0 \ | \ \rho(I + \tau \cA) < 1\, \bigr\} & \ , \ & \sigma(\cA) < - \varepsilon\, .
\end{array}
\right.
\end{equation}
\end{definition}
Thus, if $\sigma(\cA) < - \varepsilon$, then $\rho(I + \tau \cA) < 1$ for all $\tau < K(\cA, \varepsilon)$.
This value has the following meaning. If we do not know the Lyapunov exponent $\sigma (\cA)$, but have some
lower bound for $K(\cA, \varepsilon)$, then we take arbitrary $\tau < K(\cA, \varepsilon)$ and compute
the joint spectral radius $\rho(I + \tau \cA)$. If it is bigger than or equal to one, then $\sigma \ge - \varepsilon$;
otherwise, as we know, $\sigma < 0$. We have
$$
K(r, \varepsilon) \ = \ \inf\, \Bigl\{ \ K(\cA, \varepsilon) \  \Bigl| \   \max_{A \in \cA}\rho(A)\, \le \, r  \ \Bigr\}\, .
$$
Note that for some families $\cA$,
the individual bounds for $K(\cA, \varepsilon)$ can be much better than those provided by Theorems~\ref{th5} and~\ref{th10}. We are going to see this in numerical examples in Section~6. We derive lower bounds for $K(\cA, \varepsilon)$ by
analysing the Lyapunov norm of the family~$\cA$ (Proposition~\ref{p10} and~\ref{p30}), which leads to optimization problem~(\ref{prob1}) on exponential polynomials. The value of this problem is estimated in terms of Markov-Bernstein inequalities for exponents
(Proposition~\ref{p20}). Then it remains to find the infimum of those lower bounds over all families~$\cA$, this is
done in Theorem~\ref{th100}.
\smallskip

Let $\bh = (h_1, \ldots , h_d)$ be a vector such that $0 < h_1 \le h_2 \le \cdots \le h_d$ and
$\cP_{\, \bh}$ is the corresponding space of exponential polynomials $p(t) = \sum_{k  = 1}^{d}p_ke^{-h_k t}$.
For a given~$\varepsilon > 0$, we denote by~$\kappa\, (\bh, \varepsilon)$ the value of the following
minimization problem:
 \begin{equation}\label{prob1}
 \left\{
 \begin{array}{l}
 \frac{1 - p(0)}{p\, '(0) - \varepsilon p(0)} \ \to \ \min , \\
 \mbox{subject to}: \\
 p \in \cP_{\, h}\, , \\
  \|p\|_{\re_+} \le 1\, , \ p\, '(0) \, > \, \varepsilon p(0)\, .
 \end{array}
 \right.
 \end{equation}
The geometric meaning of this problem  is clarified in Proposition~\ref{p10} below.
To formulate it we need some more notation.
For a given matrix~$B$ and for $x \in \re^d$, we denote $G_B(x) =
{\rm co_s}\, \{e^{\, t\, B}x\, , \ t \in [0, +\infty)\}$, where ${\rm co_s} (X) \, = \,
{\rm co}\, \{X, -X\}$ is the symmetrized convex hull of $X.$ Thus, $G_B(x)$ is the convex hull of the curve
$\{\gamma (t)\, = \, e^{\, t\, B}x\, , \ t \in [0, +\infty)\, \bigr\}\, $ and of its reflection through the origin. If the matrix~$B$ is Hurwitz, i.e., the real parts of all its eigenvalues are negative, then
the set~$G_B(x)$ is bounded, and the curve~$\gamma$ connects the point~$x = \gamma (0)$ with the origin~$0 = \gamma(+\infty)$.

\begin{proposition}\label{p10}
For every matrix $B$ with a real negative spectrum, for every $x \in \re^d$ and $\varepsilon > 0$, the following holds:
the largest number $\tau$ such that $x + \tau \, (B\, - \, \varepsilon I)\, x \, \in \, G_{B}(x)$ is equal to
$\kappa\, (\bh, \varepsilon)$, where $\, \bh \, = \, -{\rm sp}(B)$.
\end{proposition}
{\tt Proof.} We assume without loss of generality that $B$ has $d$ distinct eigenvalues (the assertion for general matrices follows
 by taking the limit). By the Caratheorory theorem, a point belongs to the convex set $G_{B}(x)$ if
and only if this point is a convex combination of at most $d+1$ extreme points of that set.
Each extreme point of the set $G_{B}(x)$ has the form $\pm \, e^{\, t\, Bx}\, , \ t \ge 0$.
Hence, there are $n \le d+1$ nonnegative numbers $\{t_k\}_{k=1}^n$ and $n$ numbers $\{q_k\}_{k=1}^n$
such that
\begin{equation}\label{eq10}
x \ + \ \tau \, (B\, - \, \varepsilon I)\, x \quad = \quad \sum_{k=1}^n \, q_k \, e^{\, t_k\, B}\, x\, , \qquad
\sum_{k=1}^n |q_k| = 1\, .
\end{equation}
 Now let us pass to the
basis of eigenvectors of the matrix~$B$. In this basis we denote $B \, = \, {\rm diag} \, (-\beta_1, \ldots , -\beta_d)\, ,\, \beta_i > 0\, , $
$ \ e^{\, t\, B} \, = \, {\rm diag} \, (e^{\, -\beta_1 t}, \ldots , e^{-\beta_d t})$,
and $x = (x_1, \ldots , x_d)$. We can assume that all coordinates of $x$ are nonzero;
the assertion for general $x$ will again follow by taking the limit. Writing~(\ref{eq10}) coordinatewise,
we obtain
$$
x_j \, \bigl(\, 1 \, - \, \tau \, (\beta_j \, + \, \varepsilon\, )\, \bigr)  \ = \  \sum_{k=1}^n \, q_k \, e^{\, - t_k\, \beta_j}x_j,
$$
or, eliminating $x_j$:
$$
 1 \, - \, \tau \, (\beta_j \, + \, \varepsilon\, )  \quad = \quad  \sum_{k=1}^n \, q_k \, e^{\, - t_k\, \beta_j}\, .
$$
This equality does not involve~$x$. Thus, the largest~$\tau$ such that
$x + \tau (B\, -\, \varepsilon I) x \, \in \, G_{B}(x)$ is the same for all $x \ne 0$. Taking $x = \be$
(the vector of ones) we observe
that for every $\delta > 0$ the assertion $\be + \delta\, (B \, - \, \varepsilon I)\, \be \, \notin \, G_{B}(\be)$
is equivalent to the existence of a linear functional $\bep = (p_1, \ldots , p_d) \in \re^d$ separating the point $\be + \delta (B \, - \, \varepsilon I)\, \be$ from the set~$G_{B}(\be)$, i.e.,
$$
\bigl(\bep\, , \, \be \, + \, \delta\, (B \, - \, \varepsilon I)\, \be\,  \bigr) \quad > \quad \max_{y \in G_{B}(\be)} \bigl(\bep\, , \, y \bigr)\, .
$$
This follows from the convex separation theorem. The right hand side is equal to
$$
\sup_{t \in \re_+} \ \bigl|\, \sum_{k=1}^d \, p_k e^{\, - t \, \beta_k }\, \bigr| \ = \ \|p\|_{\re_+},
 $$
where $\, p(t) \, = \, \sum_{k=1}^d\, p_k e^{\, - t \beta_j }$
is an exponential polynomial. On the other hand,
$$
\bigl(\, \bar p\, , \, \be  +  \delta (B \, - \, \varepsilon I)\, \be\, \bigr)\  =\
\sum_{k=1}^d p_k \, + \, \delta \, \sum_{k=1}^d p_k(-\beta_k - \varepsilon)\  =
\  p(0) \, + \, \delta \,  p\, '(0) - \delta \, \varepsilon p(0)\, .
 $$
Normalizing, we get $\|p\|_{\re_+} \le 1$ and $ p(0) + \delta \,  p\, '(0) - \delta \varepsilon p(0)> 1$. Since $p(0) \le  \|p\|_{\re_+} = 1$, we conclude that $p\, '(0) > \varepsilon p(0)$, and hence $\delta  >  \frac{1 - p(0)}{p\, '(0) - \varepsilon p(0)}$. Thus, the point $\be \, + \, \delta \, (B \, - \, \varepsilon I)\, \be$ does not belong to
 $G_{B}(\be)$ if and only if~$\, \delta \,  > \, \frac{1 - p(0)}{p\, '(0) - \varepsilon p(0)}$, which completes the proof.

{\hfill $\Box$}
\medskip

\begin{remark}\label{r10}
{\em In fact we have proved a bit more: for every~$\tau < \kappa(\bh, \varepsilon)$ and $x \ne 0$,
the point $x + \tau \, (B\, - \, \varepsilon I)\, x$ is in the interior of~$G_B(x)$.}
\end{remark}

The value $K(\cA,  \varepsilon)$ is estimated from below by the values of~$\kappa(\bh, \varepsilon)$
with $\bh$ replaced by the spectra of the matrices $A \in \cA$ shifted by $\varepsilon$. This is done in the following proposition by analysing the Lyapunov norm of the family~$\cA$.
 \begin{proposition}\label{p30}
If all matrices from~$\cA$ have real spectra, then
\begin{equation}\label{sr}
\, K(\cA, \varepsilon) \quad \ge \quad  \min_{A \in \cA}\, \kappa\, (\bh, \varepsilon)\, , \
\end{equation}
$\mbox{where}\  \bh \ = \ -{\rm sp}(A) \, - \, \varepsilon\, \be\, .$
 \end{proposition}
{\tt Proof.} We need to show that if $\sigma (\cA) < -\varepsilon$, then $\rho
 (I + \tau A) < 1$, for all $\tau$ smaller than the right-hand side of~(\ref{sr}).
 It is well known~\cite{O, MP, B} that $\sigma (\cA) < -\varepsilon$ implies that there exists a
norm in~$\re^d$ ({\em Lyapunov norm}) such that $\|x(t)\| \, \le \, e^{-\varepsilon t} \|x_0\|$, for every trajectory~$x(\cdot)$.
In particular, this holds for a trajectory without switching, i.e., for a constant control function.
Thus, for every $A \in \cA$ we have $\|e^{\, t A}x_0\| \, < \, e^{-\varepsilon t} \|x_0\|$, and hence
$\|e^{\, t\, (A + \varepsilon I)}x_0\|\, <  \, \|x_0\|$ for all $t \in \re_+$. Therefore,
for any point $y$ from the symmetrized convex hull of the set $\, \bigl\{\, e^{\, t\, (A + \varepsilon I)}x_0 \ \bigl| \ t \in \re_+\, \bigr\}$,
we have $\|y\| \le \|x_0\|$. Applying now Proposition~\ref{p10}  for the matrix~$B = A + \varepsilon I$, and taking into account Remark~\ref{r10},
we see that
 $\, (I \, + \, \tau A) x_0$ is an interior point of the set $G_{A + \varepsilon I}(x_0)$, and, consequently,
$\, \|(I \, + \, \tau A) x_0\| \, < \,  \|x_0\|$. This means that the norm of the operator $I \, + \, \tau A$
is smaller than~$1$, for each~$A \in \cA$. Whence, $\rho (I + \tau \cA) <  1$.

{\hfill $\Box$}
\medskip

\begin{remark}\label{r20}
{\em The value~$\kappa (\bh, \varepsilon)$
of problem~(\ref{prob1}) is inversely proportional to its parameters: for every $\lambda > 0$ we have
\begin{equation}\label{homog}
\kappa \bigl(\lambda \bh\, , \, \lambda \varepsilon  \bigr) \ = \ \lambda^{-1}\, \kappa \bigl(\bh, \varepsilon \bigr)\, .
\end{equation}
To see this it suffices to change variable $u = \lambda t$ and note that ${p_t}\,'(0) = \lambda \, {p_u}'(0)$.
Hence, in the problem of computing or estimating of~$\kappa (\bh, \varepsilon)$, we can always assume that
$\bh \in \Delta$, i.e., $h_d \le 1$, otherwise, we set $\lambda = h_d^{-1}$ and apply~(\ref{homog}).
}
\end{remark}

\smallskip

Our next step is to evaluate~$\kappa(\bh, \varepsilon)$, i.e., to solve problem~(\ref{prob1}).
Let us recall (see Section~3) that  $\nu_2$ denotes the smallest positive point of alternance  of the $\bh$-Chebyshev polynomial~$T = T_{\bh}$. Thus, $T(\nu_2) = 1$,
and~$T(\cdot)$ is increasing and concave on the segment~$[0, \nu_2]$. For an arbitrary $\varepsilon >0$ we consider the following problem:
 \begin{equation}\label{prob3}
 \left\{
 \begin{array}{l}
 \frac{1 - T(t)}{T\, '(t) - \varepsilon T(t)} \ \to \ \min , \\
 \mbox{subject to}: \\
t \in [0, \nu_2]\, , \ T\, '(t) \, > \, \varepsilon T(t)\, .
 \end{array}
 \right.
 \end{equation}
In contrast to problem~(\ref{prob1}), which minimizes the functional over a unit ball in a $d$-dimensional
 space $\cP_{\bh}$, this problem is univariate and is easily solvable just by finding a
unique root of the derivative of the objective rational function $\frac{1 - T(t)}{T\, '(t) - \varepsilon T(t)}$
on the segment~$[0, \nu_2]$. Provided, of course, that the
Chebyshev polynomial $T$ is available.

\begin{proposition}\label{p20}
For every  $\varepsilon > 0$ and $\, \bh \in \Delta$
the value of problem~(\ref{prob3}) is equal to $\kappa (\bh, \varepsilon)$. Moreover,
 $$
 \kappa (\bh, \varepsilon) \quad > \quad \frac{2\, \varepsilon}{M_2(\bh)\, +\, 2\, \varepsilon^2}
 $$
\end{proposition}
The proof is in Appendix. Thus, the values of problems~(\ref{prob1}) and~(\ref{prob3}) coincide.
 This allows us to compute $\kappa (\bh, \varepsilon)$ just by evaluating
 the Chebyshev polynomial $T = T_{\bh}$ numerically and solving the simple extremal problem~(\ref{prob3})
 for it. This is done in the next section. Moreover, $\kappa (\bh, \varepsilon)$ can be estimated from below
by merely estimating the value~$T''(0)$. Our uniform lower bounds for all vectors~$\bh \in \Delta$,
are based on the following Theorem~\ref{th100} which is, in a sense, analogous to Theorem~\ref{th3} in Section~3. It states that
the smallest value is achieved when all $h_i$ take the biggest possible value~$1$, i.e.,
for $\bh = \be = (1, \ldots , 1)$.
\smallskip

\begin{theorem}\label{th100}
For a fixed $\varepsilon > 0$, the smallest value of $\kappa (\bh, \varepsilon)$
over all $\bh \in \Delta$ is attained
for $\bh = \be$, at a unique optimal polynomial $p = T_{\be} = S_d$.
\end{theorem}
The proof in in Appendix. Now we are ready to prove Theorem~\ref{th5}.

\smallskip

{\tt Proof  of Theorem~\ref{th5}.}
Combining Proposition~\ref{p10} and~\ref{p30}, we conclude that
$K(\cA, \varepsilon)$ is larger than or equal to the minimal value of $\kappa(\bh, \varepsilon)$
for $\bh = -{\rm sp}(A) - \varepsilon \be$ over all~$A \in \cA$. Theorem~\ref{th100} implies that this
value is minimal when $h_1 = \ldots = h_d$, in which case
$\bh \, = \, (\rho(A) -\varepsilon)\, \be$. Changing variables $t' = (\rho (A) - \varepsilon)t$
and invoking  Proposition~\ref{p20},  we see that
$$
\kappa(\bh, \varepsilon) \ \ge \ \frac{2\, \varepsilon}{(\rho(A) - \varepsilon)M_{2, d} \, +\, \frac{2\varepsilon^2}{\rho(A) - \varepsilon}}\, .
$$
 If $\rho(A) \le r$ for all $A \in \cA$, then we finally get
$$
K(r, \varepsilon) \ \ge \ \frac{2\, \varepsilon}{(r - \varepsilon)^2M_{2, d}  \, +\, 2\, \varepsilon^2}
$$
Since $M_{2, d}  > 2$ and $r \ge \varepsilon$, the denominator of this fraction is smaller than
$r^{\, 2} M_{2, d} $, which concludes the proof.

{\hfill $\Box$}
\medskip

\begin{remark}\label{r30}
{\em In Conjecture~\ref{conj10} we suppose that our main results hold for general matrices, and the real spectra assumption can be omitted. Proving this, probably, requires a different technique. The fact is, our approach uses essentially that any finite collection of real exponents constitute a Chebyshev system, which is
not the case for complex exponents. Actually, Propositions~\ref{p10} and \ref{p30} are true for complex
$h_k = \alpha_k + i \beta_k$, in which case we consider the space of polynomials generated by functions
$e^{-\alpha_k t}\sin \beta_k t\, , \, e^{-\alpha_k t}\cos \beta_k t$. The proofs for this case are the same.
The difficulties emerge in the solution of problem~(\ref{prob1}). In the general case we
do not know the optimal polynomial. Besides, the extremality result in Theorem~\ref{th100}, if true in general
  (which we believe), has to be proved in a completely different way, not relying on properties of Chebyshev systems.
}
\end{remark}
\bigskip

\begin{center}
\textbf{6. The algorithm and numerical examples}
\end{center}
\medskip

\begin{algorithm}\label{algo-continuous}
\DontPrintSemicolon
\KwData{A set of matrices $\cA,$ each with real spectrum, a desired accuracy $\epsilon>0.$}
\KwResult{Outputs a numerical value $s$ such that $s-\epsilon \leq \sigma(\cA) \leq s+\epsilon.$}
\Begin{
\nl Set $a,b$ lower and upper bounds on $\sigma(\cA)$\;
\nl \emph{\% for instance: $a=a_0:=\max{\mbox{Re}{(\lambda)}:\lambda \mbox{ is an eigenvalue of some } A\in \cA}$,}\;
\nl \emph{\% $b=b_0:=\max{||A||: A\in \cA}.$}\;
\nl \While{$b-a>2 \epsilon$}{
\nl set $c:=(a+b)/2$\;
\nl \label{line-tau} Take $\tau<\min\{k(A-cI,\epsilon/6): A\in\cA\}$\;
\nl \emph{\% For instance, using bounds from Theorem 2, or Theorem 3, or using the more accurate Proposition 2 if the Chebyshev polynomial is available.}\;
\nl set $\delta:= \tau \epsilon/6$\;
\nl \label{line-delta} compute an approximation $\rho^*$ of $\rho(I+\tau(A-cI-(\epsilon/6)I))$ with an accuracy $\delta$\;
\nl \label{line-approx} \emph{\% i.e. $\rho^*-\delta<\rho<\rho^*$ (For instance using \cite[Theorem 2.12]{jungers_lncis})}\;
\nl \eIf{$\rho^*-\delta>1$}{
\nl \label{line-assert1} \emph{\% $\rho>1$ and thus $\sigma(A-cI)>0$}\;
\nl set $a:=c$\;}
{\nl \label{line-assert2} \emph{\% $\rho<1+\delta$ and thus $\sigma(A-cI)<\epsilon/3$ }\;\nl set $b:=c+\epsilon/3$ }\;}
\nl {\bf Output} $(a+b)/2$
\caption{The general algorithm}}
\end{algorithm}

\normalsize

%
%
\begin{theorem}
Given a set of $d\times d$ matrices with real spectra, and maximal spectral abscissa smaller than $r,$ Algorithm \ref{algo-continuous} returns an approximation $\sigma^*$ of $\sigma(A)$ with an accuracy $\epsilon:$
$$ \sigma -\epsilon \leq \sigma^* \leq \sigma+\epsilon.$$
Moreover, the algorithm terminates within
$$
O\, \left(\, \log{\frac{L}{\epsilon}}m\, n^{\frac{3 \log m}{\log{(1-cB^2)}}}\log{\frac{12}{c B^2}}\right),
$$
where $L=N_{max}-r_{max}$ is the difference between the maximal norm and the maximal spectral abscissa of the matrices in $\cA,$ $B=\epsilon/(dr),$ and $c$ is a constant.
\end{theorem}
{\tt Proof.} 
{\bf Correctness.}
The algorithm proceeds by bisection, keeping an upper and a lower bound on $\sigma(\cA).$  It is obvious that $N_{max}$ and $r_{max}$ are valid initial values for these bounds.  Now, it remains to prove that if the property at line \ref{line-approx} is satisfied, then the two statements on Lines \ref{line-assert1} and \ref{line-assert2} of Algorithm \ref{algo-continuous} are true.

We prove the statement in Line \ref{line-assert1} by contraposition: If $\sigma(A-cI)<0,$ then, there exists an invariant convex set which contains all the trajectories of the system, and by Proposition \ref{p10} the discretized system leaves the same convex set invariant. Thus,  $$\rho(I+\tau(A-cI-(\epsilon/6)I))<1.$$

Suppose now that $\rho(I+\tau(A-cI-(\epsilon/6)I))<1+\delta$ as in Line \ref{line-assert2}.  By homogeneity of the JSR and basic arithmetic, we obtain
$$\rho(I+\frac{\tau}{1+\delta}(A-cI-(\frac{\epsilon}{6}+\frac{\delta}{\tau})I))<1.$$ By Theorem A, this implies that $\sigma(A-cI)<\frac{\epsilon}{6}+\frac{\delta}{\tau}\leq\epsilon/3.$

{\bf Running time.}\\
Line \ref{line-delta} requires an approximation of $\rho$ with absolute error bounded by $\delta= \tau\epsilon/6.$ Since $\rho$ is close to one in the worst case, (say, $\rho<2$) it is sufficient to require a relative accuracy of $\epsilon'=\delta/2=\tau\epsilon/12.$
For any arbitrarily small $\epsilon',$ algorithms are known, which deliver an approximation of the joint spectral radius with relative error $\epsilon'$ in a number of steps bounded by  $O(m\,n^{\frac{3 \log m}{\log{(1-\epsilon')}}}\log{1/\epsilon'}),$ (see \cite[Theorem 2.12]{jungers_lncis} and the proof of this theorem). Thus, a single run of the subroutine computing an approximation of $\rho$ takes \begin{equation}\label{eq-runningtime1}O(m\,n^{\frac{3 \log m}{\log{(1-\tau\epsilon/12)}}}\log{12/(\tau \epsilon)}).\end{equation}
 Recall from Theorem \ref{th5} that\footnote{In extenso, there exists a constant $C$ such that $k(r,\epsilon)>C\epsilon/(d^2r^2).$} $k(r,\epsilon)\approx\omega(\epsilon/(d^2r^2))$. Finally, the bisection algorithm divides by two the length of the interval $(a,b),$ and adds $\epsilon/3$ to this length in the worst case.  Thus, it is straightforward that the initial length $L=b_0-a_0$ becomes smaller than $2\epsilon$ within at most $O(\log{L/\epsilon}).$  Combining these last two inequalities with (\ref{eq-runningtime1}), one obtains the claim.\\

{\hfill $\Box$}
\begin{ex}
{\em Consider the LSS with the following family $\cA$ of  two $3\times 3$ matrices: 
\begin{eqnarray}
\cA&=&\left\{  \begin{pmatrix}

    -0.0622  &  0.0349 &  -0.1182\\
    0.0953 &  -0.0697  & -0.1719\\
    0.0787  &  0.0223 &  -0.2581\end{pmatrix},
    \right.   \\&&\left.
\begin{pmatrix}
 0.1591 &   0.1397 &  -0.0916\\
    0.0338  & -0.1569 &  -0.0707\\
    0.7417  &  0.3028 &  -0.4421
\end{pmatrix}\right \}.
\end{eqnarray}
We denote these matrices $A_0,A_1$ respectively.
The CQLF method gives an upper bound $\sigma\leq 0.02,$ while the spectral abscissa of both matrices is $-0.08$. 
 Thus, we have the following bounds for the Lyapunov exponent:  $-0.08\leq\sigma(\cA)\leq 0.02$, which still 
  leaves both opportunities for the system: to be stable or not. With the CQLF method, we cannot say more. 
  
  Now we apply our algorithm for computing the Lyapunov exponent.  
Both matrices happen to have real spectra, and hence we can apply our refine our analysis, thanks to the sharper estimate of a valid discretization time. 
 We choose the following parameters: we will decide whether $\sigma(\cA)>-0.079$ with an accuracy $ 0.055$ (this corresponds to a value of $\varepsilon = 3 \cdot 0.055=0.165$ in Algorithm \ref{algo-continuous}).  Applying the developments above, we
conclude that it is sufficient to discretize with the step equal to $\tau = 0.1457$ (like in Line \ref{line-tau} of the algorithm) and the maximal error $\delta$ for the JSR computation smaller or equal to $4.0\dots\ 10^{-3}$ (Line~\ref{line-delta}). (Figure \ref{fig-cheb} represents the Chebyshev polynomials corresponding to the spectra of $A_0$ and $A_1$ respectively.)
\begin{figure}\label{fig-cheb}
\centering
\begin{tabular}{cc}
\includegraphics[scale = 0.3]{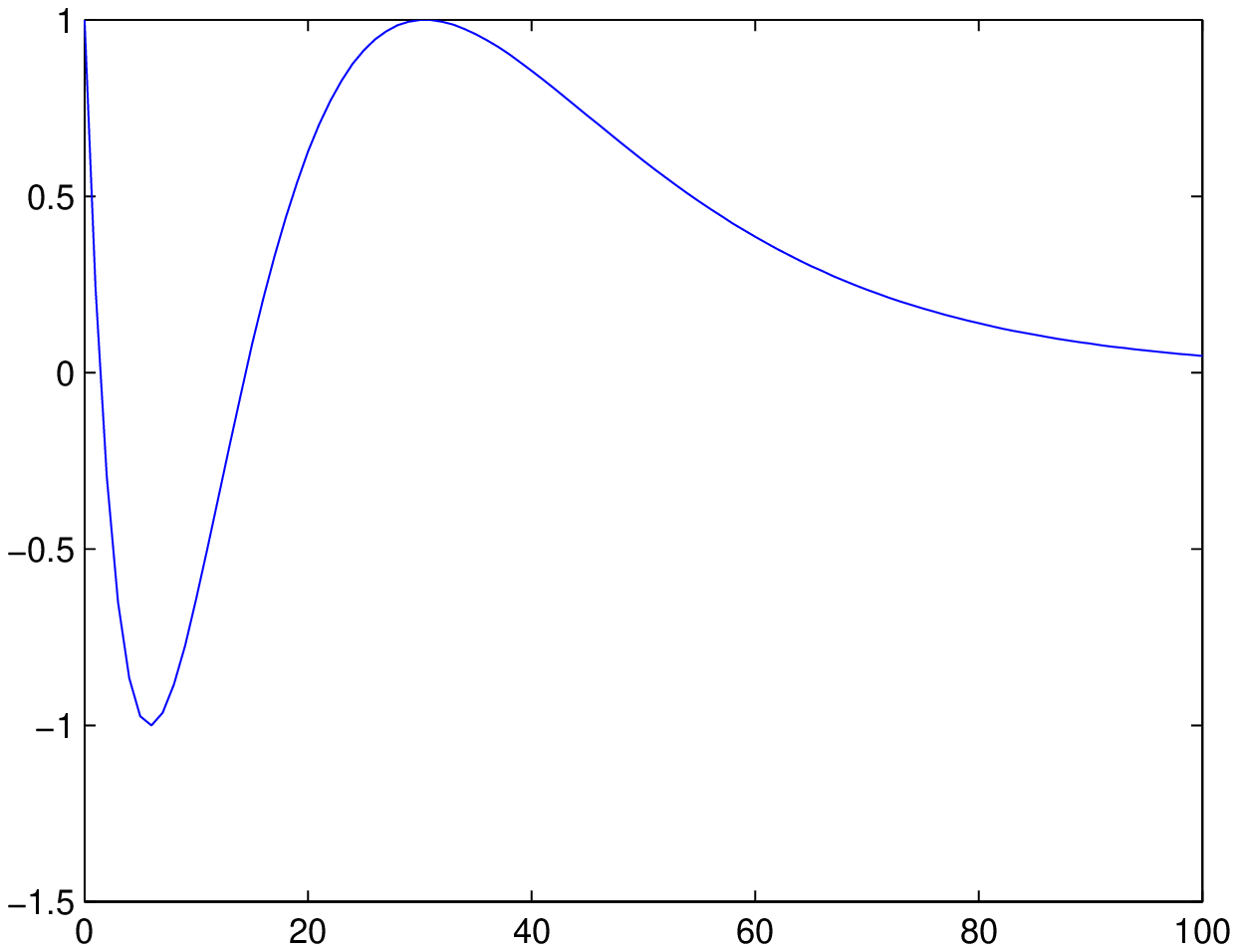}
&\includegraphics[scale = 0.3]{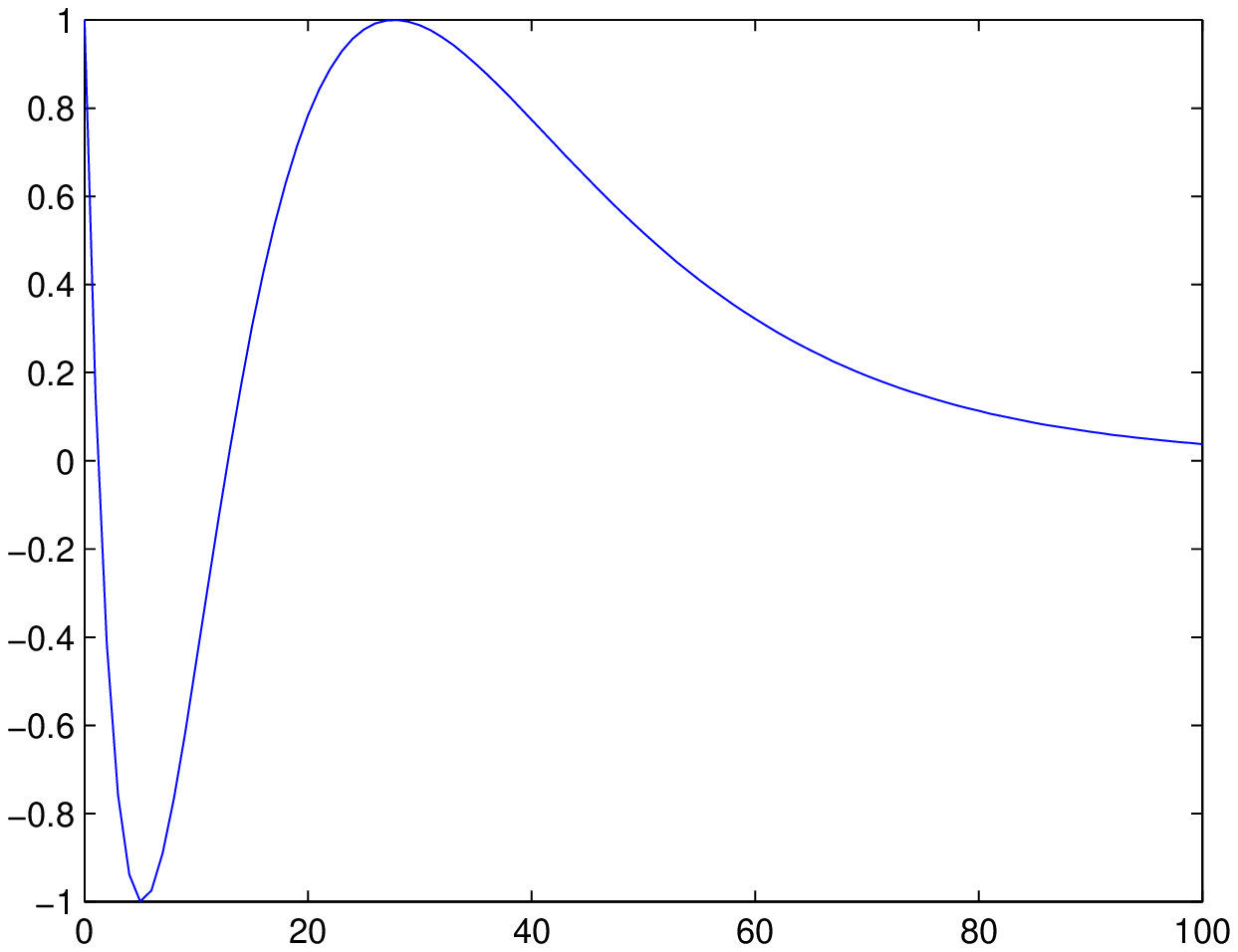}
\end{tabular}
\caption{Chebyshev polynomials corresponding to the spectra of $A_0$ and $A_1.$}
\end{figure}

%
We find that the set of matrices $I+\tau ((\cA-(-0.079\, I))-(\varepsilon/6)I)$ has the joint spectral radius smaller than $\ 1+ 3.9 \ 10^{-3}.$  (This has been obtained on a standard desktop computer with the JSR toolbox \cite{jsr-toolbox}.)  As a consequence, $\sigma(\cA-(-0.079))<0.055,$ and thus, $\sigma(\cA)<-0.024.$ This implies, in particular, that the LSS is stable. }
%
%
%
\end{ex}

\textbf{7. Conclusion}

The goal of this paper was to provide a way to compute the maximal rate of growth of a trajectory of a continuous time linear switching system, with a bound on the computation time necessary to do it with a specified accuracy.  We showed that this is possible for matrices with real spectrum, and we leave open the question for general matrices. Our techniques have a different flavor than the ones previously proposed in the literature, as they mainly aim at computing a discretization step, in order to apply efficient methods for discrete time systems (like the ones implemented in the toolbox \cite{jsr-toolbox}).  Our method is also applicable in practice for deciding stability of a LSS
and for computing the Lyapunov exponent, as demonstrated in Section 6. In the proofs we reveal a curious
link between the problem of stability of LSS and Markov-Bernstein inequalities for exponential polynomials.
As an auxiliary result, we derive a universal upper bound for the constants in those inequalities,
which is, probably, of some independent interest.

\bigskip

\textbf{Acknowledgements}. The research reported here was carried out
when the first  author was visiting the department of Mathematical Engineering,
Universit\'e Catholique de Louvain (UCL), Belgium. He
would like to thank the university for hospitality. The authors are also
grateful to Yuri Nesterov and to Konstantin Ryutin for helpful  discussions.

\bigskip

\begin{center}
\textbf{8. Appendix}
\end{center}
\medskip

{\tt Proof of Proposition~\ref{p20}.} By the compactness argument it is easily shown that problem~(\ref{prob1})
achieves its solution on some polynomial~$p \in \cP_{\bh}$. Let us prove that $p$ possess
$d-1$ points of alternance on $\re_+$: there are positive points $\mu_2 < \cdots < \mu_{d}$ such that
$p(\mu_k) = (-1)^k$. In this proof it will be more convenient to start enumeration from~$2$ than from~$1$. Take the smallest point $t_2 > 0$ such that $p(t_2) = 1$ (assume it  exists),
then the smallest point $t_3 > t_2$ such that $p(t_3) = -1$, then the
smallest point $t_4 > t_3$ such that $p(t_4) = 1$, etc. Let $t_m$ be the last point of this sequence.
If $p$ does not have the required alternance, then $m \le d-2$. Writing $t_1 = 0$, take some points $u_k \in (t_{k-1}, t_k)\, , k \ge 2,$ so that the half-intervals $[u_k, t_k)\, , \, k = 2, \ldots , m$ do not contain extremal points
of the polynomial~$p$. Set $u_1 = 0$. Since there are in total $m \le d-1$ points $u_k$, it follows that there is a polynomial
$q \in \cP_{\bh}$ such that $q(u_k) = 0, \, k = 1, \ldots , m$ and $q\, '(0)>0$. For sufficiently small number
$r> 0$ we have $\|p - rq\|_{\re_+}\le 1$. On the other hand, $(p - rq)(0) = p(0)\, $ and $\, (p - rq)\, '(0) > p\, '(0)$.
Hence the value of problem~(\ref{prob1}) is smaller for the polynomial~$p +rq$, which contradicts the assumption on~$p$.
If $t_2$ does not exists, i.e., $p(t) < 1$ for all $t >0$, then we can take $q(t) = p(0)(1-e^{-h_1t})$
and repeat the proof for this~$q$. Thus, $p$ has $d-1$ point of alternance $\mu_2 < \cdots < \mu_{d}$.
The derivative $p\, '(t)$ vanishes at all these $d-1$ points, hence it do not have other zeros on
$\re$. Therefore, $p$ increases monotone on the half-line $[-\infty, \mu_2]$. Since $|p(t)| \to \infty$
as $t \to \-\infty$, it follows that there is a unique point $\mu_1 \le 0$ such that $p(\mu_1) = -1$.
Then the polynomial $p(t -\mu_1)$ has $d+1$ points of alternance, and, by the uniqueness, it coincides
with the $\bh$-Chebyshev polynomial~$T_{\bh} = T$. Thus, $p(\cdot ) = T(\cdot+\mu_1)$.  Writing problem~(\ref{prob1}) for
$T(\cdot+\mu_1)$ and denoting $\mu_1 = t$, we arrive at~(\ref{prob3}).

It now remains to show that the value of problem~(\ref{prob3}) is larger than $\frac{2\varepsilon}{-T''(0) +2\varepsilon^2}$.
 Consider  a new function $f(t) = 1 - T(\nu_2 - t)$ on the segment~$[0, \nu_2]$ and note that $f$ is increasing and concave on $[0, \nu_2]$,
 $f(0) = f\, '(0) = 0, \|f''\|_{\re_+} \le M_2 (\bh)$. With this function,
problem~(\ref{prob3}) reads
 \begin{equation}\label{prob4}
 \left\{
 \begin{array}{l}
 \frac{f(t)}{f\, '(t) \, - \, \varepsilon \, \bigl( 1\, -\, f(t)\, \bigr)} \ \to \ \min , \\
 \mbox{subject to}: \\
t \in [0, \nu_2]\, , \ f\, '(t) \, > \, \varepsilon \, \bigl( 1\, -\, f(t)\, \bigr)\, .
 \end{array}
 \right.
 \end{equation}
Fix some $t=a$ and the value $f(t) = c$. We maximize $f\, '(t)$, in order to minimize the objective function in~(\ref{prob4}), by solving the following optimal control problem:
 $$
 \left\{
 \begin{array}{l}
 \varphi\, '(t)\ \to \ \max \\
 \mbox{subject to}\\
\varphi \in W^2_1[0, a], \ \varphi(0) \, = \, \varphi\, '(0) \, = \, 0, \, \varphi(a) = c, \\
0\, \le \, \varphi''(t) \, \le \, M_2(\bh)\ \mbox{for all}\ t \in [0, a]\, .
\end{array}
\right.
$$
The global maximum is attained for the following piecewise-quadratic function:
$\varphi(t) = 0, t \in [0, \xi]\, , \ \varphi(t) = \frac{M_2(\bh)\, t^2}{2}\, , \ t \in [\xi, a]$, where
$\xi$ is some switching point. Changing variables $t \mapsto t - \xi$, we conclude that
substituting $f(t) = \frac{M_2(\bh)\, t^2}{2}$ into~(\ref{prob4}) does not increase the value of the problem.
Now the problem becomes $\frac{\frac{M_2(\bh)t^2}{2}}{M_2(\bh) t \, - \, \varepsilon \bigl(1 -  \frac{M_2(\bh) t^2}{2}\bigr)} \, \to \min$.
Writing the determinant we see that the minimum of this function under the assumption $M_2(\bh) t \, - \, \varepsilon \bigl(1 -  \frac{M_2(\bh) t^2}{2}\bigr) \, \ge \, 0$ is equal to $\, \frac{2\, \varepsilon}{M_2(\bh) + 2\varepsilon^2}$.

{\hfill $\Box$}
\medskip

\smallskip

In the proofs of Theorems~\ref{th3} and~\ref{th100} are realized in the same way. We use the following fact, which is a
simple consequence of the convexity of norm.
\begin{lemma}\label{l20}
If $p$ and $q$ are elements of a normed space, $\|p\| = 1$ and
$\|p + \lambda_0 q\| < 1$ for some positive $\lambda_0$, then
there is a positive constant $c_0$ such that $\|p + \lambda q\|\, \le \, 1\, - \, c_0 \lambda$
for all $\lambda \in [0, \lambda_0]$.
\end{lemma}
We first give a proof of Theorem~\ref{th100} and then show how to modify it to
derive  Theorem~\ref{th3}.
\smallskip

{\tt Proof of Theorem~\ref{th100}}. Let $n$ be the biggest integer such that  $h_1 = \ldots = h_n$.
For some small $\delta > 0$ denote $\bh_{\delta} = \bigl(h_1+\delta, \ldots , h_n+\delta, h_{n+1}, \ldots , h_d\bigr)$
($\delta$ is added to the first $n$ equal entries).
Let us  show that for all sufficiently small $\delta > 0$, we have
\begin{equation}\label{repl}
\kappa\, \bigl(\bh_{\delta}, \varepsilon \bigr) \ <\
 \kappa\, \bigl(\bh, \varepsilon \bigr)\, .
\end{equation}
Thus, if $h_1 < 1$, then one can always slightly increase the smallest exponent $h_1$ to reduce~$\kappa (\bh, \varepsilon)$.
On the other hand, the set of exponential polynomials $p$ of $d$ terms such that $\|p\|_{\re_+}\le 1$, and
whose exponents are in the segment  $[h_1, 1]$  is compact. Hence,
the minimal value of $\kappa(\bh, \varepsilon)$ is attained when $h_1 =1$, i.e., at the point $\bh = \be$.
In this case a unique optimal polynomial is $T_{\be} = S_d$.

To prove~(\ref{repl}) we take a positive $\lambda$ and consider the $k$th term ($k\le n$) in the
polynomial $p(t)$. This is $p_kt^{k-1}e^{-h_1t}$. A small variation $\lambda y_k$ of the coefficient
$p_k$ and a variation $\delta \lambda$ of the exponent~$h_1$ gives
\begin{equation}\label{var1}
(p_k  +  \lambda  y_k)\, t^{k-1}e^{-\, (h_1 +  \lambda  \delta)\, t} \ - \   p_k\, t^{k-1}e^{-\, h_1t}
\ = \ \lambda \, y_k t^{k-1}e^{-h_1t}\, - \, \lambda \, \delta\, p_k t^{k}e^{- h_1t}\, + \, O(\lambda^2)\, , \quad \lambda \to 0\, .
\end{equation}
We spot the linear part of the variation using the expansion $e^{\, - \lambda \delta \, t}\, = \,
1\, - \, \lambda \delta \, t \, + \, O(\lambda^2)$. Since $p(t)$ tends to zero as $t \to \infty$,
the last term $O(\lambda^2)$ in~(\ref{var1}) is uniform over all $t \in [0, +\infty)$.
The crucial observation here is that the largest power of $t$ increases by $1$. Thus,
a small variation of the Chebyshev system $e^{-h_1t}, te^{-h_1t}, \ldots , t^{n-1}e^{-h_1t}$
leads to a larger Chebyshev system $e^{-h_1t}, te^{-h_1t}, \ldots , t^{n}e^{-h_1t}$. Now it remains to choose
coefficients $\{y_k\}_{k=1}^n$,  $\delta$ and $\lambda$ so that this variation reduce the value
of~$\kappa (\bh, \varepsilon)$.

In the proof of Proposition~\ref{p20} we showed that the extremal polynomial $p$
for problem~(\ref{prob1})   possesses $d-1$ points of alternance  $\mu_2 < \cdots < \mu_{d}$ such that
$p(\mu_k) = (-1)^k$. Every interval $[\mu_k, \mu_{k+1}]$ for $k = 2, \ldots , d-1$
 contains a unique root $u_k$ of the
polynomial~$p$. Thus, we have $d-2$ roots $u_2, \ldots , u_{d-1}$.
We add the function $t^ne^{-h_1t}$ to the system~$\cP_{\bh}$ and obtain a T-system $\cP_{\tilde \bh}$
of $d+1$ elements. Thus, we have a $(d+1)$-tuple $\tilde \bh$ such that
$\tilde h_k = h_1, \, k = 1, \ldots , n+1$, and $\, \tilde h_k = h_{k-1}, \, k = n+1, \ldots , d+1$.
Denote $u_0 = 0$ and take an arbitrary  $u_1 < \mu_2$. Since the system~$\cP_{\tilde \bh}$
contains $d+1$ element, there is a unique, up to multiplication by a constant, polynomial $q \in \cP_{\tilde \bh}$
that vanishes at the points $u_0, \ldots , u_{d-1}$ and is positive on the interval $(u_0, u_1)$.
At each point of alternance $\mu_k$ the values $p(\mu_k)$ and $q(\mu_k)$ have different signs. Therefore, for all sufficiently small $\lambda > 0$, we have $\|p +\lambda q\|_{\re_+} \, <\, 1$, and hence, by Lemma~\ref{l20},
$\|p +\lambda q\|_{\re_+} \, <\, 1\, - \, c_0\lambda$. Now let us choose coefficients
$\{y_k\}_{k=1}^n$ and   $\delta$ so that the linear part of the corresponding variation of  the polynomial
$p$ coincides with the polynomial~$q$. For $k = n+1, \ldots , n$, this is simple: we take $y_k = q_{k+1}$.
For $k = 1, \ldots , n$, we apply~(\ref{var1}), put formally $p_0 = 0$ and obtain the following system of linear equations:
$$
\begin{array}{l}
q_{1}\ = \ y_{1} \, - \, \delta\, p_{0} \\
\cdots \\
q_{k}\ = \ y_{k} \, - \, \delta\, p_{k-1} \\
\cdots \\
q_{n}\ = \ y_{n} \, - \, \delta\, p_{n-1}\\
q_{n+1} \ = \ - \, \delta\, p_{n}\, .
\end{array}
$$
Note that $p_n \ne 0$. Otherwise $p$ is a linear combination of a T-system of $d-1$ functions,
hence it cannot have an alternance of $d-1$ points. Solving the last equation:
$\delta \, = \, -\,  \frac{q_{n+1}}{p_n}$ we obtain  successively
$\, y_{k}\, = \, q_k\, - \, \frac{q_{n+1}p_{k-1}}{p_n}$. For every $\lambda > 0$
we have a polynomial $p_{\lambda}(t)\, = \,
(p_k \, + \, \lambda y_k)\, t^{k-1}e^{-(h_k+\lambda\delta)t}\, + \, \sum_{k=1}^n (p_k \, + \, \lambda y_k)\, e^{-h_kt}$
(the exponents in the last sum may also have multiplicities, in which case they
are multiplied by the corresponding powers of $t$). The linear (in $\lambda$) part of the difference $p_{\lambda}(t) \, - \,
p(t)$ coincides with $\lambda \, q(t)$. Whence, $\|p_{\lambda} - (p +\lambda q) \|_{\re_+} \, \le \, C\, \lambda^2$, and by the triangle inequality we get
$$
\|p_{\lambda}\|_{\re_+}\ \le \ \|p +\lambda q\|_{\re_+}\ +\ C\, \lambda^2\ \le \
1 - c_0\lambda \ + \ C\lambda^2 \ < \ 1
$$
 for all $\lambda> 0$ small enough. Thus,  $\|p_{\lambda}\|_{\re_+} < 1$.
Furthermore, $p_{\lambda}(0) = (p +\lambda q)(0) = p(0)$ and, since $q\, '(0) > 0$, we have
$p^{\, '}_{\lambda}(0)\, = \, p(0)\, + \, \lambda q(0) \, + \, O(\lambda^2)\, > \, p(0)$ for $\lambda> 0$ small enough.
Thus, $p_{\lambda}$ have the same value at zero, but a larger derivative. Therefore, for $p_{\lambda}$
the objective function of problem~(\ref{prob1}) is smaller than for~$p$.

It remains to show that we actually increase the exponent $h_1$, i.e., that $\delta > 0$. We have $\delta \, = \, - \,
\frac{q_{n+1}}{p_n}$. The largest term of the polynomial $p$, asymptotically as $t \to \infty$,
is $p_nt^{n-1}e^{-h_1t}$. Hence for large $t$ we have $\, {\rm sign}\, p(t)\, = \, {\rm sign}\, p_n\, = \,
(-1)^{d-1}$.  The largest at $+\infty$ term of the polynomial $q$ is
 $q_{n+1}t^{n}e^{-h_1t}$, hence for large $t$ we have ${\rm sign}\, q(t)\, = \, {\rm sign}\, q_{n+1}\, = \,
(-1)^{d}$. Thus, $p_n$ and $q_{n+1}$ have different signs, and so, $\, \delta \, = \, - \,
\frac{q_{n+1}}{p_n}\, > \, 0$.

{\hfill $\Box$}
\medskip

{\tt Proof of Theorem~\ref{th3}} is literally the same as the proof of Theorem~\ref{th100} above,
replacing the first derivative~$p'(0)$ by the $k$th derivative~$p^{(k)}(0)$. In the proof of Theorem~\ref{th100}
we showed that if $h_1 < 1$, then there exists a small perturbation $\bh_{\delta}$ of the
vector $\bh$ that reduces the norm of $p$ and increases the value $|p'(0)|$. The same perturbation
actually increase~$|p^{(k)}(0)|$ for every $k$ (the proof is the same). This implies that the maximal
$k$th derivative is achieved when $h_1 = 1$, i.e., when $\bh = \be$, which completes the proof
of Theorem~\ref{th3}.

{\hfill $\Box$}

\bigskip

\end{document}